\documentclass[reqno]{amsart}
\usepackage{amsfonts,amscd,graphicx}

\makeatletter                                                  %
\def\th@plain{\slshape}                                        %
\makeatother                                                   %

\allowdisplaybreaks

\newcommand{\tnorm}{\star}
\newcommand{\oi}{[0,1]}
\newcommand{\tot}{\leftrightarrow}
\newcommand{\llgroup}{$\ell$-group}
\newcommand{\llgroups}{$\ell$-groups}

\newcommand{\Zbb}{\mathbb{Z}}
\newcommand{\Qbb}{\mathbb{Q}}
\newcommand{\Rbb}{\mathbb{R}}
\newcommand{\Ical}{\mathcal{I}}
\newcommand{\Zcal}{\mathcal{Z}}
\newcommand{\ffrak}{\mathfrak{f}}
\newcommand{\gfrak}{\mathfrak{g}}
\newcommand{\pfrak}{\mathfrak{p}}
\newcommand{\qfrak}{\mathfrak{q}}
\newcommand{\mfrak}{\mathfrak{m}}
\newcommand{\Bfrak}{\mathfrak{B}}
\newcommand{\Luk}{\L ukasiewicz}

\newcommand{\lland}{\:\text{\&}\:}

\newcommand{\To}{\Rightarrow}

\def\dotminus{\buildrel\textstyle\cdot\over\relbar}

\newcommand{\newword}[1]{\textsl{#1}}
\newcommand{\algebra}[1]{#1}
\newcommand{\variety}[1]{\mathbf{#1}}
\newcommand{\operator}[1]{\mathbf{#1}}

\newcommand{\vect}[3]{#1_#2,\ldots ,#1_#3}

\newcommand{\abs}[1]{\lvert#1\rvert}

\newcommand{\generated}[1]{( #1 )}

\DeclareMathSymbol{\upharpoonright}{\mathrel}{AMSa}{"16}
\let\restriction\upharpoonright
\DeclareMathSymbol{\nmid}{\mathrel}{AMSb}{"2D}

\DeclareMathOperator{\Spec}{Spec}
\DeclareMathOperator{\EqTh}{EqTh}
\DeclareMathOperator{\diag}{diag}
\DeclareMathOperator{\Free}{Free}

\DeclareMathOperator{\Fl}{F\ell}
\DeclareMathOperator{\card}{card}

\theoremstyle{plain}
\newtheorem{theorem}{Theorem}[section]
\newtheorem{lemma}[theorem]{Lemma}
\newtheorem{proposition}[theorem]{Proposition}

\theoremstyle{definition}
\newtheorem{definition}[theorem]{Definition}
\newtheorem{example}[theorem]{Example}
\newtheorem{remark}[theorem]{Remark}

\begin{document}

\bibliographystyle{plain}

\sloppy

\title{Generic substitutions}

\author[G. Panti]{Giovanni Panti}
\address{Department of Mathematics\\
University of Udine\\
Via delle Scienze 208\\
33100 Udine, Italy}
\email{panti@dimi.uniud.it}

\begin{abstract}
Up to equivalence, a substitution in propositional logic is an endomorphism of its free algebra. On the dual space, this results in a continuous function, and whenever the space carries a natural measure one may ask about the stochastic properties of the action. In classical logic there is a strong dichotomy: while over finitely many propositional variables everything is trivial, the study of the continuous transformations of the Cantor space is the subject of an extensive literature, and is far from being a completed task. In many-valued logic this dichotomy disappears: already in the finite-variable case many interesting phenomena occur, and the present paper aims at displaying some of these.
\end{abstract}

\keywords{algebraic logic, substitution, stochastic properties, spectral spaces}

\thanks{\emph{2000 Math.~Subj.~Class.}: 03B50; 37B05}

\maketitle

\section{Preliminaries}

We work in the context of propositional algebraic logic.
A \newword{continuous t-norm} (t-norm for short) is a continuous function $\tnorm$ from $\oi^2$ to $\oi$ such that $(\oi,\tnorm,1)$ is a commutative monoid for which $x\le y$ implies $z\tnorm x\le z\tnorm y$. Each t-norm induces a \newword{residuum} $\to$ via $x\le (y\to z)$ iff $x\tnorm y\le z$. One checks easily that the usual lattice operations on $\oi$ are term-definable by $x\land y=x\tnorm(x\to y)$ and $x\lor y=\bigl((x\to y)\to y\bigr)\land\bigl((y\to x)\to x\bigr)$; therefore the determination of a continuous t-norm induces a uniquely defined structure of residuated lattice on $\oi$.

Let $L$ be a propositional language whose set of symbols for connectives contains $\tnorm,\to,1$. Let $M$ be a structure for $L$ satisfying:
\begin{itemize}
\item[(A1)] the support of $\algebra{M}$ is a Borel subset of $\oi$ which is either finite or of Lebesgue measure $1$;
\item[(A2)] the interpretations of $\tnorm$ and $\to$ are the restrictions to $M$ of a continuous t-norm and its residuum;
\item[(A3)] the constant $1$ is interpreted in the number $1\in M$, while the other connectives ---if any--- are interpreted as functions on $M$ of the appropriate arity that are Borel with respect to the topology that $M$ inherits from $\oi$;
\item[(A4)] if $M$ is finite, then $0\in M$ and $L$ contains a constant for $0$.
\end{itemize}
Note that the residuum is Borel, since for every $0 < a < 1$ the set $\{(x,y)\in M^2:x\to y\in[0,a)\}$ is open, while $\{(x,y)\in M^2:x\to y\in (a,1]\}$ is $F_\sigma$.

We are interested in the equational logic determined by $\algebra{M}$, i.e., the set $\EqTh(\algebra{M})$ of all term-identities $t(\bar{x})=s(\bar{x})$ that hold in $\algebra{M}$. Let $\operator{V}(\algebra{M})$ be the variety generated by $\algebra{M}$, i.e., the set of all algebras for $L$ that satisfy $\EqTh(\algebra{M})$. Observe that $\algebra{M}$ satisfies the quasi-identity
$$
x\to y=1 \lland y\to x=1  \implies x=y
$$
(indeed, assuming the premise, we have $x=x\tnorm(x\to y)=y\tnorm(y\to x)=y$). Therefore, by~\cite[Corollary~1.9]{GummUrsini84}, $\operator{V}(\algebra{M})$ is ideal-determined, and the equational logic of $\algebra{M}$ is just the algebraic counterpart of its assertional calculus~\cite{blokpigozzi89}. In the following we list a few standard examples.

\begin{example}\label{ref1}
\begin{enumerate}
\item $M$ is a chain $0=a_0 < a_1 < \cdots < a_{n-1}=1$, the t-norm $\tnorm$ is interpreted as $\min$, and a constant for $0$ is added to the language. We obtain the $n$-valued G\"odel-Dummet logic~\cite{dummett}, which is boolean logic if $n=2$;
\item $\tnorm$ and $0$ as in (1), $M=\oi$; we obtain the intuitionistic propositional calculus plus the prelinearity axiom $(x\to y)\lor(y\to x)=1$;
\item $M$ and $\tnorm$ as in (1), with a connective $\sim$ added to the language and interpreted by $\sim a_i=a_{i+1}\pmod{n}$. We get Post's $n$-valued logic~\cite{post}, \cite[\S2.3]{pantisu};
\item $M=(0,1]$, $\tnorm=$ the ordinary product of real numbers. Then the residuum gets interpreted as $a\to b=\min\{1,b/a\}$. This is similar to the Product Logic in~\cite[\S4.1]{hajek98}, but since $0\notin M$ every term induces a continuous function;
\item $M=\oi$, $\tnorm$ the \Luk\ conjunction $a\tnorm b=\max(0,a+b-1)$: this is the logic of Wajsberg hoops~\cite{AglianoPanti02};
\item as in (5), but a constant for $0$ is added to the language; we obtain the \Luk\ infinite-valued logic;
\item $M=\oi$. Denote temporarily the \Luk\ conjunction in~(5) by $\odot$, and define a t-norm $\tnorm$ on $M$ as follows:
$$
a\tnorm b=\begin{cases}
2^{-1}(2a\odot 2b), & \text{if $a,b<1/2$;} \\
2^{-1}\big(1+((2a-1)\odot(2b-1))\big), & \text{if $1/2\le a,b$;} \\
\min\{a,b\}, & \text{otherwise.}
\end{cases}
$$
The resulting logic coincides with Hajek's Basic Logic~\cite[\S2.2]{hajek98} on the set of all $1$-variable formulas~\cite{Montagna}.
\end{enumerate}
\end{example}

In all this paper, $\variety{V}=\operator{V}(M)$ will denote a variety of the above form. Every algebra $A\in\variety{V}$ carries a natural order given by $a\le b$ iff $a\to b=1$. A \newword{filter} in $\algebra{A}$ is the counterimage of $1$ under some homomorphism of domain $A$. Ideal-determinacy of $\variety{V}$ means that the lattice of congruences and the lattice of filters of $\algebra{A}$ are isomorphic via $\theta\mapsto 1/\theta$. We say that a filter $\pfrak$ of $\algebra{A}$ is \newword{prime} if, for any two filters $\ffrak,\gfrak$, if $\pfrak\supseteq \ffrak\cap\gfrak$ then $\pfrak\supseteq \ffrak$ or $\pfrak\supseteq \gfrak$. 
Since a lattice structure is term-interpretable in all algebras of $\variety{V}$, the variety $\variety{V}$ is congruence-distributive~\cite[Theorem~12.3]{burrissan81}, and $\pfrak$ is prime iff the implication
$$
\pfrak= \ffrak\cap \gfrak \implies \pfrak= \ffrak \text{ or } \pfrak= \gfrak
$$
holds.

A filter $\ffrak=\generated{a}$ is \newword{principal} if it is generated by a single element $a$; in other words, $\generated{a}$ is the intersection of all filters that contain $a$. An \newword{ideal term in $y$} is a term $t(\vect x1n,y)$ such that the equation $t(\vect x1n,1)=1$ is valid in $\variety{V}$. By~\cite[Lemma~1.2]{GummUrsini84},
$$
\generated{a}=\{t(\vect c1n,a):t(\bar{x},y)\text{ is an ideal term in $y$ and }
\vect c1n\in\algebra{A}\}.
$$

\begin{proposition}
$\generated{a}\cap\generated{b}=\generated{a\lor b}$. \label{ref2}
\end{proposition}
\begin{proof}
The right-to-left inclusion is clear. 
Let $t(\bar{x},y),s(\bar{z},w)$ be terms, not necessarily ideal. Then the equation
\begin{displaymath}
t(\bar{x},y\lor w)\land s(\bar{z},y\lor w)\le t(\bar{x},y)\lor s(\bar{z},w)
\tag{$*$}
\end{displaymath}
is true in $M$. Indeed, if $y,w$ are interpreted in $l,m$, then $y\lor w$ is interpreted either in $l$ or in $m$. Since $\variety{V}$ is generated by $M$, ($*$) holds in $\variety{V}$. Assume now that $e=t(\bar{c},a)=s(\bar{d},b)\in\generated{a}\cap\generated{b}$, for certain ideal terms $t(\bar{x},y),s(\bar{z},w)$ in $y$ and $w$, respectively, and certain $\bar{c},\bar{d}\in\algebra{A}$. Then the term $r(\bar{x},\bar{z},u)=t(\bar{x},u)\land s(\bar{z},u)$ is ideal in $u$, and we have
$$
e=t(\bar{c},a)\lor s(\bar{d},b)\ge r(\bar{c},\bar{d},a\lor b)\in\generated{a\lor b}.
$$
Since clearly filters are closed upwards, $e\in\generated{a\lor b}$ and the proof is complete.
\end{proof}

\begin{proposition}\label{ref3}
Let $\algebra{A}\in\variety{V}$. The following are equivalent:
\begin{itemize}
\item[(i)] $\algebra{A}$ is totally-ordered;
\item[(ii)] the set of filters of $\algebra{A}$ is totally-ordered;
\item[(iii)] every filter is prime;
\item[(iv)] the trivial filter $\generated{1}$ is prime.
\end{itemize}
The above conditions are implied by ---but do not imply--- subdirect irreducibility of~$\algebra{A}$.
\end{proposition}
\begin{proof}
The implications (i)$\To$(ii)$\To$(iii)$\To$(iv) are clear. Assume~(iv), and let $a,b\in\algebra{A}$. We have $1=(a\to b)\lor(b\to a)$, and therefore by Proposition~\ref{ref2} $\generated{1}=\generated{a\to b}\cap\generated{b\to a}$. Since the filter $\generated{1}$ is prime, it must be either $a\le b$ or $b\le a$, and therefore $\algebra{A}$ is totally-ordered. If $\algebra{A}$ is subdirectly irreducible, then trivially the filter $\generated{1}$ is prime. Let $\algebra{A}=(\oi,\tnorm,\to,0,1)$ be as in Example~\ref{ref1}(2), $\theta$ a nontrivial congruence on $\algebra{A}$, and let $\ffrak=1/\theta\not=\generated{1}$ be the filter associated to $\theta$. Then there exists a unique $a\in[0,1)$ such that $\ffrak$ equals either $(a,1]$ or $[a,1]$, and $\theta=\diag(\algebra{A})\cup \ffrak^2$. Indeed, assume $b\theta c$ and, without loss of generality, $b\le c$. Then $b=c\to b\in \ffrak$ and therefore both $b$ and $c$ are in $\ffrak$. By defining $a=\inf \ffrak$, our description of $\ffrak$ and $\theta$ follows. It is immediate that $\algebra{A}$ is a totally-ordered subdirectly reducible algebra.
\end{proof}

Let $\Spec A$ be the set of all proper prime filters of $A$, excluding the improper filter $A$. Endow $\Spec A$ with the hull-kernel topology, by taking the family of all sets of the form
$$
O_a=\{\pfrak\in\Spec\algebra{A}:a\notin \pfrak\}
$$
as an open subbasis. The following properties are well known~\cite{Hochster69}, \cite[Chapter~10]{bkw}, \cite{Priestley94}:
\begin{enumerate}
\item $\Spec A$ is a ---possibly noncompact--- \newword{spectral space}, i.e., a space possessing an intersection-closed basis of compact open sets, and in which every irreducible closed set is the closure of a unique point;
\item the closure of $Y\subseteq\Spec A$ is $\{\pfrak\in\Spec A:\pfrak\supseteq\bigcap Y\}$, and in particular $Y$ is dense iff $\bigcap Y=\{(1)\}$. The relation $\pfrak\le\qfrak$ iff $\pfrak\subseteq\qfrak$ iff $\qfrak$ is in the closure of $\pfrak$ is the \newword{specialization order} of $\Spec A$;
\item every filter of $A$ is a ---possibly empty--- intersection of proper prime filters (because the proper filter $\pfrak$ is prime if $A/\pfrak$ is subdirectly irreducible). Hence the mapping
$$
\ffrak\mapsto O_\ffrak=\{\pfrak\in\Spec A:\ffrak\not\subseteq\pfrak\}
$$
is an isomorphism between the lattice of filters of $A$ and the lattice of open sets in $\Spec A$;
\item under the above mapping, the finitely generated filters correspond to the compact open sets. By Proposition~\ref{ref2}, these are closed under finite intersections. In particular, $\Spec A$ is compact iff $A$ is finitely generated as a filter (an example of an algebra with a noncompact spectrum is the one used in the proof of Proposition~\ref{ref3}, as soon as one removes the constant $0$ from the language).
\end{enumerate}

\section{Substitutions acting on spectra}\label{ref8}

A \newword{syntactical substitution} is an endomorphism of the algebra of terms in $L$. Since we work up to the equational theory of $\variety{V}$, we are actually interested in endomorphisms of the free algebra in $\variety{V}$ over $\kappa\le\omega$ many generators; we denote this latter algebra by $\Free_\kappa(\variety{V})$. Endomorphisms of $\Free_\kappa(\variety{V})$ will be named \newword{substitutions over $\kappa$ variables}; a substitution is \newword{invertible} if it is an automorphism of $\Free_\kappa(\variety{V})$. It is well known that $\Free_\kappa(\variety{V})$ is the subalgebra of $\algebra{M}^{\algebra{M}^\kappa}$ generated by the projection functions $x_i:M^\kappa\to M$, for $0\le i < \kappa$. We give to $M$ the topology induced by the standard topology on $\oi$, and to $M^\kappa$ the product topology.

Let $\lambda$ be the probability measure naturally associated with $M$. More precisely, if $M$ has Lebesgue measure $1$ as a subset of $\oi$, then $\lambda$ is the restriction of the Lebesgue measure to $M$, while if $M$ is finite then $\lambda$ is the counting measure; in both cases $\lambda(M)=1$. We extend $\lambda$ to a probability measure on $M^\kappa$ by defining
$$
\lambda\bigl(C(A_{i_1},\ldots,A_{i_t})\bigr)=\lambda(A_{i_1})\times\cdots\times
\lambda(A_{i_t}),
$$
where the $A_{i_j}$s are measurable subsets of $M$, and $C(A_{i_1},\ldots,A_{i_t})$ is the cylinder $\{p\in M^\kappa:x_{i_j}(p)\in A_{i_j} \text{ for every } j=1,\ldots,t\}$.

For every $p\in\algebra{M}^\kappa$, we have the homomorphism $\Free_\kappa(\variety{V})\to\algebra{M}$ given by evaluation at $p$. Let $\pi(p)$ be the kernel, which is a prime filter by Proposition~\ref{ref3}, since $M$ is totally-ordered. For every $t\in\Free_\kappa(\variety{V})$, let $Z(t)=\{p\in M^\kappa:t(p)=1\}$ be the \newword{$1$-set} of $t$. The mapping $\pi$ is Borel: indeed, the counterimage of the basic closed set $F_t=\Spec\Free_\kappa(\variety{V})\setminus O_t$ is $Z(t)$, which is Borel by our assumption (A3). Note that $\pi$ is undefined at the points $p\in\bigcap\{Z(t):t\in\Free_\kappa(\variety{V})\}=Z_0$, since then $\pi(p)$ is the improper filter. This gives no trouble, since $Z_0$ is either empty (e.g., if $L$ contains a constant for $0$), or contains the element $p_0=(\ldots,1,\ldots)$ only, and in this latter case $\lambda(Z_0)=0$ by our assumption (A4).

We push forward $\lambda$ to a Borel probability measure on $\Spec\Free_\kappa(\variety{V})$ by setting
$$
\lambda(A)=\lambda(\pi^{-1}[A]).
$$
We use $\lambda$ both for the original measure on $M^\kappa$ and for the induced measure on $\Spec\Free_\kappa(\variety{V})$; this should cause no confusion.

We briefly recall a few basic facts of ergodic theory~\cite{Walters82}, \cite{CornfeldFomSi82}, \cite{PollicottYur98}: we state them at a level of generality appropriate to our setting. A \newword{measure-theoretic dynamical system} is a triple $(X,\mu,S)$ where $X$ is a second countable topological space, $\mu$ is a Borel probability measure on $X$ (i.e., $\mu(X)=1$ and $\mu$ is defined on the Borel $\sigma$-algebra~$\Bfrak$), and $S:X\to X$ is a measurable mapping. We let $A,B,\ldots$ vary over the elements of $\Bfrak$, and we drop universal quantifications over these objects in our definitions. We say that $S$ is \newword{nonsingular} if $\mu(A)=0$ implies $\mu(S^{-1}[A])=0$, and 
\newword{measure-preserving} if $\mu(A)=\mu(S^{-1}[A])$.
Assume that $S$ is measure-preserving. Then $S$ is:
\begin{enumerate}
\item \newword{ergodic} if $S^{-1}[A]=A$ implies $\mu(A)=0$ or $\mu(A)=1$;
\item \newword{mixing} if
$$
\lim_{n\to\infty}\mu(A\cap S^{-n}[B])=\mu(A)\cdot\mu(B);
$$
\item \newword{exact} if the tail $\sigma$-field
$\bigcap_{n\ge0}S^{-n}\Bfrak$ contains only sets of measure $0$ or $1$.
\end{enumerate}
It is a classical fact that exactness $\To$ mixing $\To$ ergodicity. We single out in the following proposition the property of ergodic systems that mainly concerns us.

\begin{proposition}\label{ref7}
Let $(X,\mu,S)$ be an ergodic system, and assume that $\mu$ is supported on all of $X$ (i.e., every nonempty open set has measure $>0$). Then for $\mu$-all points $x\in X$ the following holds:
\begin{itemize}
\item[$(*)$] for every $k\ge0$, the point $S^k(x)$ has a dense orbit.
\end{itemize}
\end{proposition}
\begin{proof}
One of several characterizations of ergodicity is that $\mu(A)>0$ implies $\mu(\bigcup\{S^{-n}[A]:0\le n\})=1$ \cite[Theorem~1.5]{Walters82}. Since $S$ is measure-preserving, for all $k\ge 0$ the set $B(A,k)=\bigcup\{S^{-n}[A]:k\le n\}$ has also measure $1$. Let $\{A_i:i<\omega\}$ be a countable basis of nonempty open sets. Then the set of all $x$ that satisfy property $(*)$ is $\bigcap\{B(A_i,k):i,k<\omega\}$, which has measure~$1$.
\end{proof}

Let $\sigma:\Free_\kappa(\variety{V})\to\Free_\kappa(\variety{V})$ be a substitution. We always assume that $\sigma$ is \newword{nontrivial}, i.e., that the counterimage of every proper prime filter is proper; this condition is automatically satisfied if $0$ is a constant of the language. Let $X$ be $\Spec\Free_\kappa(\variety{V})$, $\lambda$ the probability measure on $X$ defined at the beginning of this section, $S$ the transformation given by $S(\pfrak)=\sigma^{-1}[\pfrak]$. 
By our assumption $\kappa\le\omega$, the space $X$ is second countable.
Since $S^{-1}[F_t]=F_{\sigma(t)}$, we have that $S$ is continuous and $(X,\lambda,S)$ is a measure-theoretic dynamical system. 
Let us say that $\sigma$ is \newword{minimal} if in the associated dynamical system $(X,S)$ all points $x\in X$ have a dense orbit~\cite[Definition~5.1]{Walters82}. The minimality property (which is purely topological, since no measure is involved) is quite relevant to logic.
Indeed, let $t(\bar{x})=s(\bar{x})$ be an $n$-variable identity in the language $L$. Then $t(\bar{x})=s(\bar{x})$ is true in $\variety{V}$ iff it is true in the generators of $\Free_n(\variety{V})$ iff $t$ and $s$ are equal modulo $\pfrak$, for $\pfrak$ ranging on a dense subset of $\Spec\Free_n(\variety{V})$. Let $\sigma$ be a minimal substitution over $n$ variables. Then in order to test the truth of $t(\bar{x})=s(\bar{x})$ we may choose any prime filter $\pfrak$ and test the identity modulo 
all points of the $S$-orbit of $\pfrak$ (dually stated, we test all identities $\sigma^n(t)=\sigma^n(s)$, for $n\ge0$, modulo a fixed $\pfrak$).
For all $\pfrak$, this procedure provides a correct test.

Most many-valued logics do not admit minimal substitutions. For example, we shall see that both for \Luk\ and for product logic the space of maximal filters is homeomorphic to a closed Euclidean ball on which $S$ acts continuously; such an action necessarily has fixed points by the Brouwer fixed point theorem. We are therefore lead to taking into consideration measure-theoretic issues. 

\begin{definition}
Let\label{ref20} $X=\Spec\Free_\kappa(\variety{V}),\lambda$ be as above. We say that a Borel probability measure $\mu$ on $X$ is \newword{algebraically equivalent} to $\lambda$ if there exists an automorphism $\rho$ of $\Free_\kappa(\variety{V})$ such that, writing $R$ for the dual homeomorphism, we have $\mu(A)=\lambda(R^{-1}[A])$ for every $A$. Let $\sigma$ be a substitution, $S$ its dual. We say that $\sigma$ is \newword{generic} if for every $\mu$ algebraically equivalent to $\lambda$ the set $\{x\in X:x\text{ has a dense $S$-orbit}\}$ has $\mu$-measure $1$.
\end{definition}

\begin{remark}
The reader should compare the notion of algebraic equivalence of measures with that of topological equivalence in~\cite[\S~3]{oxtobyulam41}. The two notions coincide if all homeomorphisms of $X$ come from automorphisms of $\Free_\kappa(\variety{V})$: this is the case, e.g., of classical logic.
For the two cases that mainly concern us, namely \Luk\ logic and falsum-free product logic over finitely many variables, algebraic equivalence gives essentially no new measures. Indeed, by~\cite[Theorem~2.6]{dinolagripa} and our Theorem~\ref{ref9}, $\lambda$ is invariant under the action of the automorphism group of \Luk\ logic. Analogously, we will see in Section~4 that formulas of falsum-free product logic correspond to piecewise-linear homogeneous functions over the negative cone of $\Rbb^n$. Automorphisms of such structures are induced by finite families of nonsingular $n\times n$ matrices, as explained in~\cite{beynon77}. This fact and Theorem~\ref{ref17} imply that all measures algebraically equivalent to $\lambda$ share the same $\sigma$-ideal of nullsets. I thank Andrew Glass and Daniele Mundici for clarifying discussions over this point.
\end{remark}

We think of genericity as a replacement for minimality, for logics that do not admit minimal substitutions. When testing an identity along the orbit of a prime filter, the initial choice of the filter becomes relevant; however, from the measure-theoretic point of view, all choices are good.
Note that in certain pathological situations the condition~($*$) of Proposition~\ref{ref7} results to be stronger that genericity, since it neglects the transient behaviour of the orbits; see Example~\ref{ref5}.

Let $s_i=\sigma(x_i)$. Then the $\kappa$-tuple $(\ldots,s_i,\ldots)$ determines a function $\bar{s}:M^\kappa\to M^\kappa$ via $p\mapsto(\ldots,s_i(p),\ldots)$, and the following diagram commutes
$$
\begin{CD}
M^\kappa  @>{\bar{s}}>>  M^\kappa  \\
@V{\pi}VV                @VV{\pi}V   \\
X         @>>S>        X
\end{CD}
$$
Note that the diagram is valid even when $\pi$ is undefined at $p_0=(\ldots,1,\ldots)$. Indeed, since $\sigma$ is nontrivial we have $\bar{s}^{-1}[\{p_0\}]=\{p_0\}$, and we may safely substitute $M^\kappa$ in the top row with $M^\kappa\setminus\{p_0\}$.
The mapping $S$ is continuous, while $\bar{s}$ and $\pi$ are in general only Borel. Each of the maps $\sigma$ and $\bar{s}$ determines the other; moreover, $\bar{s}$ determines $S$ since $\pi$ has dense range. If $\pi$ is injective then $S$ determines $\bar{s}$.
This may fail if $\pi$ not injective: one easily constructs an example by using the logic of Example~\ref{ref1}(4). Its $1$-generated free algebra is the chain $1>x>x\tnorm x>x\tnorm x\tnorm x>\cdots$, whose spectrum is the single point $\generated{1}$. The substitution $x\mapsto 1$ is trivial, while every other substitution induces the identity map on the spectrum. Note however that $\pi$ is injective for many relevant cases, e.g., for classical logic (where it is even a homeomorphism) and for \Luk\ logic. Let us present here a pathological example.

\begin{example}\label{ref5}
Let $M$ be as in Example~\ref{ref1}(2), $\variety{V}=\operator{V}(M)$. Then $\Free_1(\variety{V})$ is isomorphic to the direct product $\{0,1\}\times\{0,1/2,1\}$ (the factors being subalgebras of $M$), with free generator $(0,1/2)$. This is easily seen by observing that:
\begin{enumerate}
\item for every $0<a<1$, the subalgebra of $M$ generated by $a$ is isomorphic to $\{0,1/2,1\}$;
\item for every term $t(x)$ we have $t(1)=1$ iff $t(1/2)\ge 1/2$.
\end{enumerate}
Hence the subalgebra of $\oi^{\oi}$ generated by the identity function is isomorphic to $\{0,1\}\times\{0,1/2,1\}\times\{0,1\}$ (by (1)), but the third factor is superfluous (by (2)). By the way, an analogous argument shows that $\Free_n(\variety{V})$ is a finite product of finite chains; this fact was proved in~\cite{horn69}. Under the specialization order, $X=\Spec\Free_1(\variety{V})$ is the poset
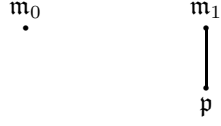
\begin{center}
\unitlength0.8cm\thinlines
\begin{picture}(3,3)(0,0)
\put(3,1){\circle*{0.1}}
\put(3,0.7){\makebox(0,0){$\pfrak$}}
\put(0,2){\circle*{0.1}}
\put(3,2){\circle*{0.1}}
\put(0,2.3){\makebox(0,0){$\mfrak_0$}}
\put(3,2.3){\makebox(0,0){$\mfrak_1$}}
\put(3,2){\line(0,-1){1}}
\end{picture}
\end{center}
where
\begin{itemize}
\item $\mfrak_0=\{1\}\times\{0,1/2,1\}=\pi(0)$;
\item $\pfrak=\{0,1\}\times\{1\}=\pi(a)$, for every $a\in(0,1)$;
\item $\mfrak_1=\{0,1\}\times\{1/2,1\}=\pi(1)$.
\end{itemize}
Note that $\lambda$ is concentrated in $\pfrak$: $\lambda(A)=1$ iff $\pfrak\in A$. Let $\sigma$ be the substitution that maps $x$ to $s=x\to 0$. Then $\bar{s}(0)=1$ and $\bar{s}(a)=0$ for $a\in(0,1]$; also
$S(\pfrak)=S(\mfrak_1)=\mfrak_0$ and $S(\mfrak_0)=\mfrak_1$. Both $\bar{s}$ and $S$ are singular. The $S$-orbit of $\pfrak$ is all of $X$, and hence $\sigma$ is generic.
On the other hand, the condition~($*$) of Proposition~\ref{ref7} is not satisfied, because the transient behaviour is relevant: the orbit of $S(\pfrak)$ is not dense.
\end{example}

Recall that a measure $\nu$ is \newword{absolutely continuous} w.r.t.\ a measure $\mu$ on the same space if $\mu(A)=0$ implies $\nu(A)=0$; we write then $\nu\ll\mu$.

\begin{definition}
Let $\sigma$ be a substitution, $S$ its associated continuous mapping on $X=\Spec\Free_\kappa(\variety{V})$. We say that $\sigma$ is \newword{ergodic} (respectively, \newword{mixing} or \newword{exact}) if there exists a Borel probability
measure $\mu$ on $X$ such that:
\begin{enumerate}
\item all measures algebraically equivalent to $\lambda$ are $\ll\mu$;
\item $(X,\mu,S)$ is ergodic (respectively, mixing or exact).
\end{enumerate}
\end{definition}

\begin{remark}
\begin{enumerate}
\item Every ergodic substitution is generic, but the converse is false: see Theorem~\ref{ref16}(3).
\item The introduction of the auxiliary measure $\mu$ is unavoidable; in Example~\ref{ref14} we shall see two ergodic substitutions whose dual maps do not preserve $\lambda$.
\item The reader may wonder why we did not define $\sigma$ to be generic if the set of points of $M^\kappa$ that have a dense orbit under $\bar{s}$ has measure $1$ under all measures algebraically equivalent to $\lambda$. This definition is simpler, since it avoids the introduction of the machinery of spectra, and stronger, since a substitution which is generic according to it is also generic in our sense. The trouble is that it is too strong: the mapping $\pi$ may collapse many points, so that $\bar{s}$ does not act generically on the points of $M^\kappa$, although it does on the fibers of $\pi$. Stated otherwise: the orbit of a random point of $M^\kappa$ may be not dense in $M^\kappa$, but still contain sufficiently many points to ascertain logical truth. We shall see an example in Theorem~\ref{ref16}(3.2).
\end{enumerate}
\end{remark}

In the remainder of this section, we discuss the familiar case of classical logic. With the obvious modifications, our discussion applies to Post's $n$-valued logic as well, because both logics are functionally complete. Let then $\variety{V}$ be the variety of boolean algebras. The mapping $\pi$ is a homeomorphism that identifies the discrete set $2^n$ and the Cantor space $2^\omega$ with $\Spec\Free_n(\variety{V})$ and $\Spec\Free_\omega(\variety{V})$, respectively. In the first case $\lambda$ is the counting measure (which is obviously fixed by algebraic equivalence), and in the second is the measure determined by
$$
\lambda\bigl(C(\vect e{{i_1}}{{i_t}})\bigr)=2^{-t},
$$
where $\vect i1t < \omega$, $\vect e{{i_1}}{{i_t}}\in\{0,1\}$, and $C(\vect e{{i_1}}{{i_t}})$ is the clopen cylinder $\{p\in 2^\omega:p(i_j)=e_{i_j}\}$. Every substitution is nontrivial, and every continuous mapping from $2^\kappa$ to itself is induced by a unique substitution. For substitutions over finitely many variables the situation trivializes.

\begin{proposition}\label{ref4}
No substitution over $n<\omega$ variables is exact or mixing. For every such substitution $\sigma$, the following are equivalent:
\begin{itemize}
\item[(i)] $\sigma$ is ergodic;
\item[(ii)] $\sigma$ is generic;
\item[(iii)] $\sigma$ is minimal;
\item[(iv)] $S$ permutes cyclically the $2^n$ elements of $\Spec\Free_n(\variety{V})$;
\item[(v)] $\sigma$ is invertible and the group it generates has order $2^n$.
\end{itemize}
\end{proposition}
\begin{proof}
If $S$ preserves a measure $\mu\gg\lambda$, then it must be surjective. Hence it is bijective and $\mu=\lambda$; clearly $S$ cannot be mixing or exact. (i) $\To$ (ii) always holds, and (ii) $\To$ (iii) holds because every point has nonzero measure. By following the orbit of a point one sees that a minimal $S$ satisfies (iv). (iv) $\To$ (v) trivially, and (iv) $\To$ (i) because $S$ is then ergodic with respect to $\lambda$. Assume (v), and decompose $S$ in the product of disjoint cycles in the symmetric group over $2^n$ letters. The order of $S$ is the least common multiple of the lengths of the cycles. Therefore each cycle has length a power of $2$, which in turn implies that there is only one cycle. Hence (iv) holds.
\end{proof}

As we noted in the Abstract, the study of continuous transformations of the Cantor space, i.e., in our language, of substitutions over $\omega$ variables in classical logic, is the subject of an extensive literature. The simplest minimal ---hence generic--- substitutions are those given by adding machines. Realize the Cantor space as the space underlying the topological group $\Zbb_p$ of $p$-adic integers, for some prime~$p$. Let $\alpha\in\Zbb_p$ be such that $\chi(\alpha)\not=1$ for every nontrivial continuous character $\chi$; this amounts to the invertibility of $\alpha$ in the local ring $\Zbb_p$. Then~\cite[pp.~97--99]{CornfeldFomSi82} guarantees that the translation by $\alpha$ ($\beta\mapsto\alpha+\beta$) is a minimal homeomorphism of $\Zbb_p$. The simplest choice is to take $\alpha=1\in\Zbb_2$, and by explicit computation one sees that the corresponding substitution $\sigma$ is given by
\begin{equation*}
\begin{split}
x_0 & \mapsto \neg x_0 \\
x_{i+1} & \mapsto \big((x_0\land x_1 \land \cdots \land x_i)\tot\neg x_{i+1}\big)
\end{split}
\end{equation*}
Note that restricting $\sigma$ to the first $n$ variables gives a substitution of the form described in Proposition~\ref{ref4}, namely the translation by $1$ in $Z_{2^n}$.

All minimal transformations of the Cantor space can be realized as mappings on the path spaces of certain combinatorial objects, called Bratteli-Vershik diagrams~\cite[Theorem~4.7]{HermanPutSk92}. Moreover, a unique dimension group with order unit~\cite{effros} can be read off either from the transformation or from a corresponding diagram~\cite[Theorem~5.4 and Corollary~6.3]{HermanPutSk92}; this dimension group is a complete invariant for strong orbit equivalence of Cantor minimal systems~\cite{GiordanoPutSk95}. It is remarkable that Bratteli diagrams and dimension groups play a r\^ole also in \Luk\ logic~\cite{mundiciam}, \cite[\S4.2]{pantites}. It is unclear whether this recurrent appearance is a coincidence or an indicator of some hidden underlying structure.

\section{\Luk\ logic}

We have seen in the previous section that for classical and Post logics the generic substitutions are either trivial (in the finite-variable case) or extremely complex (in the infinite-variable case). In this section we shall describe completely the stochastic properties of $1$-variable substitutions in \Luk\ logic. In this logic, the $1$-variable case has always been an excellent test case~\cite{Mundici92}, \cite{mundicipas}, for two main reasons:
\begin{itemize}
\item it presents the piecewise-linear structure of the free algebras in a simplified nontrivial fashion;
\item the extensions of \Luk\ logic ---i.e., the subvarieties of the variety $\variety{MV}$ of MV-algebras--- are axiomatizable by $1$-variable identities~\cite{komori}, \cite{DiNolaLettieri1999}, \cite{pantijancl}, and hence correspond bijectively to the fully invariant filters of $\Free_1(\variety{MV})$~\cite[Corollary~14.10]{burrissan81}.
\end{itemize}

In this section $\variety{MV}=\operator{V}(M)$ will denote the variety of Example~\ref{ref1}(6), whose elements are named MV-algebras. We presuppose familiarity with the basic theory of MV-algebras; see~\cite{cignolimunta}, \cite{CignoliMundici97b} for a quick introduction or~\cite{CignoliOttavianoMundici00} for a more extensive treatment.

A \newword{McNaughton function over the $n$-cube}~\cite{mcnaughton},
\cite{mundicimn} is a continuous function $t:[0,1]^n\to[0,1]$
for which the following holds:
\begin{itemize}
\item[] there exist finitely many affine linear
polynomials $\vect t1k$, each
$t_i$ of the form
$t_i=a^0_ix_0+a^1_ix_1+\cdots+a^{n-1}_ix_{n-1}+a^n_i$,
with $a^0_i,\ldots,a^n_i$ integers, such that, for each $p\in[0,1]^n$,
there exists $i\in\{1,\ldots,k\}$ with $t(p)=t_i(p)$.
\end{itemize}
It is well known that the free MV-algebra over $n$ generators $\Free_n(\variety{MV})$ is the algebra of all McNaughton functions over the $n$-cube under pointwise operations. Note that we tacitly identify an $n$-variable term with the function it induces.

Let $\sigma$ be a substitution over $n$ variables. As in Section~\ref{ref8} we have the commuting diagram
$$
\begin{CD}
\oi^n  @>{\bar{s}}>>  \oi^n  \\
@V{\pi}VV             @VV{\pi}V   \\
\Spec\Free_n(\variety{MV})    @>>S>   \Spec\Free_n(\variety{MV})
\end{CD}
$$
By~\cite[Proposition~8.1]{mundicijfa} $\pi$ maps homeomorphically the $n$-cube to the subspace of maximal elements of $\Spec\Free_n(\variety{MV})$. As an example, we draw a picture of 
$\Spec\Free_1(\variety{MV})$ under the specialization order:
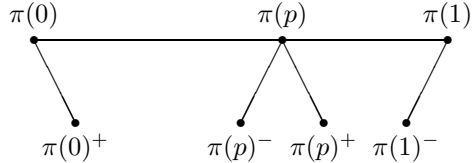
\begin{center}
\unitlength1.1cm\thinlines
\begin{picture}(5,3)(0,0)
\put(0.5,1){\circle*{0.1}}
\put(2.5,1){\circle*{0.1}}
\put(3.5,1){\circle*{0.1}}
\put(4.5,1){\circle*{0.1}}
\put(0.5,0.7){\makebox(0,0){$\pi(0)^+$}}
\put(2.5,0.7){\makebox(0,0){$\pi(p)^-$}}
\put(3.5,0.7){\makebox(0,0){$\pi(p)^+$}}
\put(4.5,0.7){\makebox(0,0){$\pi(1)^-$}}
\put(0,2){\circle*{0.1}}
\put(3,2){\circle*{0.1}}
\put(5,2){\circle*{0.1}}
\put(0,2.3){\makebox(0,0){$\pi(0)$}}
\put(3,2.3){\makebox(0,0){$\pi(p)$}}
\put(5,2.3){\makebox(0,0){$\pi(1)$}}
\put(0,2){\line(1,0){5}}
\put(0,2){\line(1,-2){0.5}}
\put(3,2){\line(1,-2){0.5}}
\put(3,2){\line(-1,-2){0.5}}
\put(5,2){\line(-1,-2){0.5}}
\end{picture}
\end{center}
For every rational point $p\in(0,1)$, there are two prime filters $\pi(p)^+$ and $\pi(p)^-$ attached below the maximal $\pi(p)$. If $p$ is irrational, then $\pi(p)$ is a minimal prime filter. Finally, there are  two points $\pi(0)^+$ and $\pi(1)^-$ attached below $\pi(0)$ 
and $\pi(1)$. The prime filter $\pi(p)^-$ is the filter of all McNaughton functions that are $1$ in a left neighborhood of $p$, and analogously for $\pi(p)^+$, $\pi(0)^+$, and $\pi(1)^-$. 
See~\cite{pantiprime} for a generalization of this description to $\Spec\Free_n(\variety{MV})$.

In the cases of classical and Post logics, we saw that  $\pi$ provides a topological conjugacy between $\bar{s}$ and $S$. This is definitely false in \Luk\ logic, since the spaces $\oi^n$ and $\Spec\Free_n(\variety{MV})$ are not even homeomorphic. Nevertheless, as we shall see in Theorem~\ref{ref9}, \Luk\ logic enjoys the nice property that $\bar{s}$ and $S$ are measure-theoretically isomorphic (provided that $\bar{s}$ is nonsingular, which is of course the interesting case). As a consequence, we can carry out our quest of stochastic properties at the simpler level of transformations of the $n$-cube.

Recall~\cite[Definition~2.4]{Walters82} that two measure-theoretic dynamical systems $(X_i,\mu_i,S_i)$, for $i=1,2$, are \newword{isomorphic} if there are Borel subsets $A_i\subseteq X_i$ of full measure such that
\begin{enumerate}
\item $S_i[A_i]\subseteq A_i$;
\item there exists a measure-preserving Borel isomorphism $\rho:A_1\to A_2$ such that $\rho\circ S_1=S_2\circ\rho$.
\end{enumerate}

\begin{theorem}
Consider\label{ref9} the systems $(\oi^n,\lambda,\bar{s})$ and $(\Spec\Free_n(\variety{MV}), \lambda, S)$ as above. Then $\bar{s}$ is nonsingular iff $S$ is nonsingular. If this happens, the two systems are isomorphic under $\pi$.
\end{theorem}
\begin{proof}
If $S$ is singular, then $\bar{s}$ is singular by the definition of the measure $\lambda$ on $\Spec\Free_n(\variety{MV})$. Consider now $\bar{s}$. By~\cite[p.~2]{mcnaughton} and \cite[Theorem~2.6]{dinolagripa} there exists a finite partition $\{P_1,\ldots,P_r\}$ of the $n$-cube in compact convex polyhedrons of  dimension $n$ such that $\bar{s}$ is expressible by
$$
\begin{pmatrix}
\alpha_0 \\
\vdots \\
\alpha_{n-1}
\end{pmatrix}
\mapsto
U_k
\begin{pmatrix}
\alpha_0 \\
\vdots \\
\alpha_{n-1}
\end{pmatrix}
+
V_k
$$
on each polyhedron $P_k$; in the above expression $U_k$ is an $n\times n$ matrix and $V_k$ a column vector, both having integer entries.
Suppose that $\bar{s}$ is singular. Then $U_k$ must be singular for some $k$, and hence $\bar{s}[P_k]$ is a polyhedron of dimension $<n$. Let $t\in\Free_n(\variety{MV})$ be such that the $1$-set $Z(t)$ of $t$ is exactly $\bar{s}[P_k]$; such a $t$ exists by the theory of Schauder hats~\cite{mundicimn}, \cite{pantilu}. Let $F_t=\Spec\Free_n(\variety{MV})\setminus O_t$ be the basic closed set determined by $t$. Then $\lambda(F_t)=\lambda(\pi^{-1}[F_t])=\lambda(Z(t))=0$. On the other hand
\begin{equation*}
\begin{split}
\lambda(S^{-1}[F_t]) &= \lambda(F_{\sigma(t)}) \\
&= \lambda(Z(\sigma(t))) \\
&\ge \lambda(P_k) \\
&> 0,
\end{split}
\end{equation*}
since $\sigma(t)=t\circ\bar{s}$ and $P_k\subseteq Z(t\circ\bar{s})$; hence $S$ is singular.

We now assume that $\bar{s}$ and $S$ are nonsingular, and prove our second statement. For every $t\in\Free_n(\variety{MV})$, let
$$
D_t=\bigcap\{U:U\text{ is an open set in }\Spec\Free_n(\variety{MV})
\text{ and }U\supseteq F_t\}.
$$

\noindent\emph{Claim.} $D_t$ is Borel (actually $G_\delta$) and $\pfrak\in D_t$ iff the maximal filter to which $\pfrak$ specializes (i.e., the unique maximal $\mfrak$ such that $\mfrak\supseteq\pfrak$) is in $F_t$.

\noindent\emph{Proof of Claim.} Let $U$ be an open set containing $F_t$, and write $U=\bigcup\{O_{r_i}:i\in I\}$. Since $F_t$ is closed and $\Spec\Free_n(\variety{MV})$ is compact (because $\Free_n(\variety{MV})$ is finitely generated as an improper filter), there exists a finite subset $J$ of $I$ such that, writing $r=\bigwedge_{i\in J}r_i$, we have $F_t\subseteq O_r\subseteq U$. This implies that $D_t$ is a (necessarily countable) intersection of basic open sets. Assume now that the maximal specialization $\mfrak$ of $\pfrak$ is not in $F_t$, and let $\mfrak=\pi(p)$ for a uniquely determined $p\in\oi^n$. Let $E$ be a closed $n$-dimensional ball centered in $p$ and having empty intersection with $Z(t)$. By~\cite[Corollary~3.4]{mundiciFree} there exists $q\in\Free_n(\variety{MV})$ such that $q\restriction E=1$ and $Z(q)\cap Z(t)=\emptyset$.
The latter identity implies $F_q\cap F_t=\emptyset$, and the former implies $q\in\pfrak$ by~\cite[Proposition~3.1]{mundiciFree}. Therefore $\pfrak\notin O_q\supseteq F_t$, and $\pfrak\notin D_t$. Since open sets are downwards closed in the specialization order, the reverse inclusion $\mfrak\in F_t$ $\To$ $\pfrak\in D_t$ is immediate.\qed

Let $\Zcal$ be the set of all $t\in\Free_n(\variety{MV})$ such that $\lambda(Z(t))=0$. Then, for every $t\in\Zcal$,
\begin{equation*}
\begin{split}
\lambda(D_t) &= \lambda\bigl(\pi^{-1}\bigl[\bigcap\{U:U
\text{ is open and }U\supset F_t\}\bigr]\bigr) \\
&= \lambda\bigl(\bigcap\{\pi^{-1}[U]:U
\text{ is open and }U\supset F_t\}\bigr) \\
&= \lambda\bigl(\bigcap\{V\subseteq\oi^n:V
\text{ is open and }V\supseteq Z(t)\}\bigr) \\
&= 0,
\end{split}
\end{equation*}
since the Lebesgue measure on $\oi^n$ is regular. Let
\begin{equation*}
\begin{split}
B &= \bigcup\{Z(t):t\in\Zcal\} \\
C &= \bigcup\{D_t:t\in\Zcal\}
\end{split}
\end{equation*}
Then $C$ is the set of all $\pfrak\in\Spec\Free_n(\variety{MV})$ whose maximal specialization $\mfrak=\pi(p)$ is such that $p\in B$. In particular $\pi^{-1}[C]=B$. Let $A_1=\oi^n\setminus B$, $A_2=\Spec\Free_n(\variety{MV})\setminus C$. Then:
\begin{itemize}
\item $p\in A_1$ iff $\pi(p)\in A_2$;
\item $A_1$ and $A_2$ are Borel sets of full measure;
\item $\pi\restriction A_1$ is surjective on $A_2$. Indeed, let $\pfrak\in A_2$ and let $\pi(p)$ be its maximal specialization. Then $p\notin B$. Write $p=(\vect\alpha0{{n-1}})\in\oi^n$. The real numbers $\vect\alpha0{{n-1}},1$ must be linearly independent over $\Qbb$ since, otherwise, we would be able to find $t\in\Zcal$ such that $t(p)=1$. By~\cite[Theorem~4.8 and Corollary~4.9]{pantiprime}, no element of $\Spec\Free_n(\variety{MV})$ specializes to $\pi(p)$ except $\pi(p)$ itself. Hence $\pfrak=\pi(p)$;
\item $\bar{s}[A_1]\subseteq A_1$ and $S[A_2]\subseteq A_2$. Indeed, we only need to show that $\bar{s}(p)\in B$ implies $p\in B$. Let $\bar{s}(p)\in Z(t)$, for some $t\in\Zcal$. Since $Z(t\circ\bar{s})=\bar{s}^{-1}[Z(t)]$ and $\bar{s}$ is nonsingular, $t\circ\bar{s}\in\Zcal$. Therefore $p\in Z(t\circ\bar{s})\subseteq B$.
\end{itemize}
Since $\pi$ is injective and measure-preserving, its restriction to $A_1$ is an isomorphism as required.
\end{proof}

In our next theorem we will characterize the generic substitutions over $1$ variable. Let $s\in\Free_1(\variety{MV})$ and let $Q_0\subset\Qbb$ be the finite set of points of nondifferentiability of $s$, including $0$ and $1$. For $i\ge 0$, write $Q_{i+1}=Q_i\cup s[Q_i]$. The chain $Q_0\subseteq Q_1\subseteq Q_2\subseteq\cdots$ stabilizes after finitely many steps. Indeed, for every rational number $p=a/b\in\oi$ in lowest terms, let us say that $b$ is the \newword{denominator} of $p$. Then, since the linear pieces of $s$ have all integer coefficients, the denominators of the elements of $Q_i$ are bounded by the maximum denominator $d$ of the elements of $Q_0$. It follows that $\bigcup_iQ_i$ is contained in the finite set of all rational numbers in $\oi$ whose denominator is $\le d$, and our claim follows.

Display the set of points of $\bigcup_iQ_i$ as $0=q_0<q_1<\cdots <q_r=1$, and let $I_i=[q_{i-1},q_i]$, for $1\le i\le r$, be the corresponding intervals. On each $I_i$ the function $s$ is of the form $a_ix+b_i$, with $a_i,b_i\in\Zbb$, and $s[I_i]=I_q\cup I_{q+1}\cup\cdots\cup I_{q+t}$, for some $1\le q\le q+t\le 1$.
Define the \newword{Markov graph} of $s$ to be the directed graph $G_s$ whose set of vertices is $\{\vect I1r\}$, and there is an edge connecting $I_i$ with $I_j$ iff $s[I_i]\supseteq I_j$. The \newword{Markov matrix} of $s$ is the $r\times r$ matrix $E_s$ whose $(i,j)$th entry is $0$ if $s[I_i]\not\supseteq I_j$ and $\abs{a_i}^{-1}$ otherwise.

A directed graph $G$ is \newword{strongly connected} if every two entries are connected by a path. The \newword{period} of a strongly connected graph is the g.c.d.\ of the lengths of the paths starting from some given vertex and returning to it; the period does not depend on the choice of the initial vertex. $G$ is \newword{primitive} if it is strongly connected and has period $1$.
By~\cite[Theorem~4.5.8]{LindMar95} $G_s$ is primitive iff there exists $k\ge 1$ such that $E_s^k$ has all entries $>0$.

\begin{example}
Consider\label{ref14} the function $s$ having graph
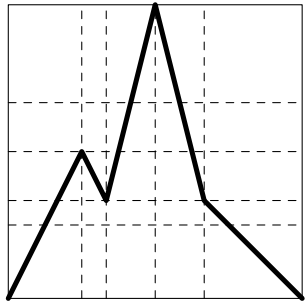
\begin{figure}[h]
\begin{center}
\includegraphics[height=4cm,width=4cm]{figura.1}
\end{center}
\end{figure}

\noindent The set of points of nondifferentiability is $Q_0=\{0,1/4,1/3,1/2,2/3,1\}$, which is stable under $s$. We have
$$
E_s=
\begin{pmatrix}
1/2 & 1/2 & 1/2 & 0 & 0 \\
0 & 0 & 1/2 & 0 & 0 \\
0 & 0 & 1/4 & 1/4 & 1/4 \\
0 & 0 & 1/4 & 1/4 & 1/4 \\
1 & 1 & 0 & 0 & 0
\end{pmatrix}
$$
and $G_s$ is primitive.

Let now $t$ have graph
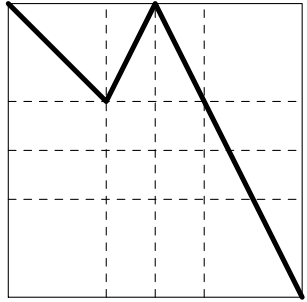
\begin{figure}[h]
\begin{center}
\includegraphics[height=4cm,width=4cm]{figura.2}
\end{center}
\end{figure}

\noindent We have $Q_0=\{0,1/3,1/2,1\}$ and $Q_1=\{0,1/3,1/2,2/3,1\}$, which is stable. The Markov matrix is
$$
E_t=
\begin{pmatrix}
0 & 0 & 0 & 1 \\
0 & 0 & 0 & 1/2 \\
0 & 0 & 0 & 1/2 \\
1/2 & 1/2 & 1/2 & 0
\end{pmatrix}
$$
and $G_t$ is strongly connected of period $2$.
\end{example}

\begin{theorem}
Let\label{ref10} $\sigma:x_0\mapsto s\in\Free_1(\variety{MV})$ be a substitution over $1$ variable in \Luk\ logic. Then the following are equivalent:
\begin{itemize}
\item[(i)] $\sigma$ is generic;
\item[(ii)] $\sigma$ is ergodic;
\item[(iii)] $\sigma$ is nonivertible and $G_s$ is strongly connected.
\end{itemize}
Furthermore, the following are equivalent:
\begin{itemize}
\item[(iv)] $\sigma$ is mixing;
\item[(v)] $\sigma$ is exact;
\item[(vi)] $\sigma$ is nonivertible and $G_s$ is primitive.
\end{itemize}
If any of (i)--(vi) holds, then the measure $\mu\gg\lambda$ with respect to which $s$ is ergodic is unique and $\mu\ll\lambda$ holds as well. The density function $d\mu/d\lambda\in L_1(\oi,\lambda)$ has rational values and is constant on each interval $\vect I1r$.
\end{theorem}

The rest of this section is devoted to the proof of Theorem~\ref{ref10}. First of all, it is not difficult to show that the only invertible substitutions over $1$ variable are the identity and the flip $x_0\mapsto (x_0\to 0)$~\cite[Example~2.7]{dinolagripa}; neither of them is generic. If no path in $G_s$ connects $I_i$ with $I_j$, then no point of $I_i$ will ever enter the topological interior of $I_j$, and $\sigma$ is not generic. This proves (i) $\To$ (iii); we already know (ii) $\To$ (i) and (v) $\To$ (iv).

We prove (iv) $\To$ (vi). Assume that $\sigma$ is mixing (and hence generic), and let $t$ be the period of the strongly connected graph $G_s$. By~\cite[\S4.5]{LindMar95} the vertices of $G_s$ can be partitioned in $t$ equivalence classes $\vect\Ical1t$ such that the quotient graph (defined by $\Ical_i\to\Ical_j$ iff there exists an edge from some element of $\Ical_i$ to some element of $\Ical_j$) has a cyclic structure $\Ical_1\to\Ical_2\to\cdots\to\Ical_t\to\Ical_1$. By using the fact that $s$ maps the extreme points of the real unit interval to themselves, one sees easily that $t\le2$. Assume by contradiction that $t=2$, and let $A_i=\bigcup\Ical_i$. Then $s^{-1}[A_1]=A_2$ and $s^{-1}[A_2]=A_1$. This clearly implies that $s$ is not mixing with respect to any measure $\gg\lambda$, which is a contradiction. We conclude that $t=1$ and (vi) holds.

\begin{lemma}
Assume\label{ref11} (iii) in Theorem~\ref{ref10}. Then there exists a 
density function $f\in L_1(\oi,\lambda)$ such that:
\begin{itemize}
\item $f$ takes strictly positive rational values, and is constant on each interval $\vect I1r$;
\item the measure $\mu$ determined by $d\mu=f\,d\lambda$ is $s$-invariant.
\end{itemize}
\end{lemma}
\begin{proof}
Every nonsingular $s$ has finite fibers, and determines the \newword{Perron-Frobenius operator}~\cite[Chapter~3]{LasotaMac94}
$$
P:L_1(\oi,\lambda)\to L_1(\oi,\lambda)
$$
via
$$
(Pf)(x)=\sum_{y\in s^{-1}[x]}\frac{f(y)}{\abs{s'(y)}}.
$$
If $f\in L_1(\oi,\lambda)$ has constant value $a_i$ on each $I_i$, then we identify $f$ with the row vector $(a_1\,\cdots\,a_r)$, and $Pf$ with the row vector $(a_1\,\cdots\,a_r)E_s$. By~\cite[Theorem~1.3.5]{Kitchens98}, $E_s$ has a real eigenvalue $\alpha>0$ and a corresponding left eigenvector $(a_1\cdots a_r)$ (unique up to scalar multiples) such that:
\begin{itemize}
\item $\alpha\ge\abs{\beta}$, for every eigenvalue $\beta$;
\item $a_i>0$, for every $i$.
\end{itemize}
Write $(b_1\cdots b_r)^{tr}$ for the column vector defined by $b_i=\lambda(I_i)$, and normalize $(a_1\cdots a_r)$ by setting $\sum a_ib_i=1$. Let $f\in L_1(\oi,\lambda)$ be the density function corresponding to $(a_1\cdots a_r)$. Since $s$ is surjective, we have
$$
\alpha=\begin{pmatrix}
a_1 & \cdots & a_r
\end{pmatrix}
E_s
\begin{pmatrix}
b_1 \\
\vdots \\
b_r
\end{pmatrix}=
\int_{s^{-1}\oi}(Pf)\,d\lambda=
\int_{\oi}f\,d\lambda=1,
$$
and $f$ is a fixed point for $P$ (in particular, $(a_1\cdots a_r)\in\Qbb^r$). As a consequence, the measure $\mu$ determined by
$$
\mu(A)=\int_Af\,d\lambda
$$
is $s$-invariant. Since $f$ is never $0$, each of $\mu$ and $\lambda$ is absolutely continuous with respect to the other.
\end{proof}

Note that if a transformation is ergodic w.r.t.~a measure $\nu$, then it does not preserve any other measure $\ll\nu$. In particular, if $s$ is ergodic w.r.t.~$\mu\gg\lambda$, then necessarily $\mu$ is the measure given by Lemma~\ref{ref11}.
The invariant densities for the functions in Example~\ref{ref14} are those corresponding to the vectors $(6/5,6/5,9/5,3/5,3/5)$ and $(3/4,3/4,3/4,3/2)$, respectively.
We present the following well known theorem (see~\cite[p.~290]{CornfeldFomSi82}, \cite[Theorem~V.2.2]{deMelovan93}, \cite[Theorem~12.5]{PollicottYur98}) in a simplified form which is convenient for our needs.

\begin{theorem}
Let\label{ref22} $h$ be a surjective map form $\oi$ to itself, let $0=q_0<q_1<\cdots<q_r=1$, and let $I_i=[q_{i-1},q_i]$ for $1\le i\le r$. Assume that:
\begin{enumerate}
\item $h$ is $C^2$ on $U=\oi\setminus\{\vect q0r\}$ and there exist $\alpha,\beta>0$ such that $\alpha\le\abs{h'(x)}$ and $\abs{h'(x)},\abs{h''(x)}\le\beta$ for every $x\in U$;
\item there exist $k\ge1$ and $\gamma>1$ such that $\abs{(h^k)'(x)}\ge\gamma$ for every $x$ in which the derivative is defined;
\item every $h[I_i]$ is the union of ---necessarily consecutive--- intervals in $\{\vect I1r\}$, and the resulting Markov graph is primitive.
\end{enumerate}
Then we have:
\begin{itemize}
\item[(a)] there exists a unique $h$-invariant probability measure 
$\mu\ll\lambda$. Its density $d\mu/d\lambda$ is strictly positive and uniformly bounded away from zero; in particular $\lambda\ll\mu$ as well;
\item[(b)] $h$ is exact with respect to $\mu$;
\item[(c)] $\mu(A)=\lim_{n\to\infty}\lambda(h^{-n}[A])$ for every Borel $A\subseteq\oi$.
\end{itemize}
\end{theorem}
\begin{proof}
In order to apply~\cite[Theorem~V.2.2]{deMelovan93}, we only have to check that the conditions on~\cite[p.~353]{deMelovan93} are satisfied. This boils down to showing that:
\begin{itemize}
\item[(4)] there exist $\varepsilon,C>0$ such that every $h\restriction I_i$ is a $C^{1+\varepsilon}$ diffeomorphism satisfying
$$
\biggl\lvert\frac{h'(x)}{h'(y)}-1\biggr\rvert\le
C\cdot\abs{h(x)-h(y)}^\varepsilon;
$$
\item[(5)] there exists $K>0$ and $\delta>1$ such that
$$
\abs{(h^n)'(x)}\ge K\delta^n,
$$
for every $n\ge1$ and every $x$ for which the derivative is defined.
\end{itemize}
The bounds in (1) assure that every $h\restriction I_i$ is a $C^{1+\varepsilon}$ diffeomorphism, for $\varepsilon=1$. We obtain~(4) from the intermediate value theorem: $\abs{h'(x)/h'(y)-1}=\abs{h'(y)}^{-1}\cdot\abs{h'(x)-h'(y)}\le\alpha^{-1}\abs{h'(x)-h'(y)}\le\alpha^{-1}\beta\abs{x-y}\le\alpha^{-2}\beta
\abs{h(x)-h(y)}$. If we can take $\alpha>1$ in~(1), then~(5) is immediate. Otherwise, let $K=\gamma^{-1}\alpha^{k-1}$ and $\delta=\gamma^{1/k}$. Writing $n=uk+v$ for $0\le v<k$ we obtain
\begin{equation*}
\begin{split}
(h^n)'(x) &= (h^v)'(x)\cdot
\prod_{j=0}^{u-1}(h^k)'(h^{jk+v}(x)) \\
& \ge \alpha^{k-1}\gamma^u=K\gamma^{u+1}\ge K\delta^n.
\end{split}
\end{equation*}
\end{proof}

Assume now (vi) in Theorem~\ref{ref10}. We claim that $\abs{(s^2)'(x)}\ge2$ whenever the derivative is defined. 
For every interval $I_i$, let $0\not= a_i\in\Zbb$ be the value of the derivative of $s$ in $I_i$. No $a_i$ can be $1$, for otherwise the graph $G_s$ would not be strongly connected. Let $I_i\to I_j$ be an arrow in $G_s$; we just need to show that $\abs{a_ia_j}\ge2$. If this was not the case, then $a_i=a_j=-1$ and $s\restriction I_i=s\restriction I_j=-x+1$. But then one sees easily that $G_s$ contains the arrow $I_j\to I_i$ as well, and no other arrow starting from $I_i$ or from $I_j$. This is impossible since $G_s$ is primitive, and our claim is settled. We can then apply Theorem~\ref{ref22}, thus obtaining (vi) $\To$ (v). To complete the proof of Theorem~\ref{ref10}, we assume that (iii) holds and prove that $s$ is ergodic w.r.t.\ the measure $\mu$ in Lemma~\ref{ref11}. If $G_s$ is primitive we are done. Otherwise, by the discussion following the statement of Theorem~\ref{ref10}, the intervals $\vect I1r$ are partitioned in two equivalence classes $\Ical_1,\Ical_2$ such that the quotient graph is cyclic. Consider $s^2$; for every edge $a\to b$ in $G_s$, let $I_{a,b}=I_a\cap s^{-1}[I_b]$. Then the intervals $I_{a,b}$ are the basic intervals for $s^2$. Moreover, $I_{a,b}\to I_{c,d}$ is an edge in $G_{s^2}$ iff $s^2[I_{a,b}]\supseteq I_{c,d}$ iff $s[I_b]\supseteq I_{c,d}$ iff $s[I_b]\supseteq I_c$ iff $I_b\to I_c$ is an edge in $G_s$. This means that $G_{s^2}$ is the disjoint union of two graphs $G_1$ and $G_2$. The vertices of $G_1$ are the edges of $G_s$ whose starting vertex is in $\Ical_1$, and two such vertices $a\to b$, $c\to d$ are connected by an edge in $G_1$ iff $b\to c$ is an edge in $G_s$; a dual description holds for $G_2$.

\begin{lemma}
Both\label{ref13} $G_1$ and $G_2$ are primitive.
\end{lemma}
\begin{proof}
We prove the statement for $G_1$. Let $a\to b$, $c\to d$ be vertices of $G_1$. Since $G_s$ is strongly connected, there is a path in $G_s$ connecting $b$ with $c$. By looking at the edges of this path alternatively as edges and vertices of $G_1$, we obtain a path in $G_1$ connecting $a\to b$ with $c\to d$; hence $G_1$ is strongly connected. All vertices of a strongly connected graph $G$ have the same period; this means that the period of $G$ can be defined as the g.c.d.\ of the lengths of the simple closed circuits in $G$. In our case, the circuits of $G_1$ correspond bijectively to the circuits of $G_s$, and this correspondence doubles the lengths. $G_s$ has period $2$, and hence $G_1$ must have period $1$.
\end{proof}

The intervals in $\Ical_1$ (respectively, $\Ical_2$) must be consecutive; indeed, suppose by contradiction $I,J\in\Ical_1$, and $K\in\Ical_2$ is between $I$ and $J$. For some odd $k\ge1$ we have $s^k[K]\supseteq I\cup J$ and then, since $s^k$ is continuous, $s^k[K]\supseteq K$, which is impossible. Let $\bigcup\Ical_1=[0,q]$ and $\bigcup\Ical_2=[q,1]$. Again by Lemma~\ref{ref13} and Theorem~\ref{ref22}, $s^2\restriction[0,q]$ and $s^2\restriction[q,1]$ are both exact with respect to the appropriate restrictions of the measure $\mu$ in Lemma~\ref{ref11}. We prove that $s$ is ergodic w.r.t.\ $\mu$ by showing that $\lambda(A)>0$ implies $\lambda(\bigcup\{s^{-n}A:s\ge0\})=1$. Without loss of generality, $B=A\cap[0,q]$ has nonzero measure, and hence so does $s^{-1}B\subseteq[q,1]$. We have
$$
\bigcup_{n\ge0}s^{-n}A \supseteq
\bigcup_{n\ge0}s^{-n}B =
\bigcup_{k\ge0}s^{-2k}B \cup
\bigcup_{k\ge0}s^{-2k-1}B.
$$
The two sets to the right intersect at most in $q$, and have full measure in $[0,q]$ and in $[q,1]$, respectively. Therefore, $\bigcup\{s^{-n}A:n\ge0\}$ has full measure in $\oi$. This completes the proof of Theorem~\ref{ref10}.

\section{Falsum-free product logic}

In this final section we will discuss the product logic of Example~\ref{ref1}(4). We will show that the spectrum of its $n$-generated free algebra is homeomorphic to $\Spec\Free_{n-1}(\variety{MV})$, with preservation of the measure $\lambda$. The substitutions over $n$ variables give rise to continuous piecewise-fractional transformations, and the resulting dynamics is  richer than the one in \Luk\ logic. In Theorem~\ref{ref16} we will see that such substitutions may have attracting fixed points, or may be generic without being ergodic.

Recall from Example~\ref{ref1}(4) that the falsum-free product logic is defined by the structure $M=((0,1],\tnorm,\to,1)$, where $\tnorm$ is the product of real numbers and $a\to b=\min\{1,b/a\}$. In all this section $\variety{V}$ will denote the variety generated by~$M$.
We assume that the reader is familiar with the basic theory of lattice-ordered abelian groups (\llgroups) \cite{fuchs63}, \cite{birkhoff}, \cite{bkw}, \cite{darnel95}, in particular with the description of the free \llgroup\ over $n$ generators $\Fl(n)$ in terms of continuous piecewise-linear homogeneous functions with integer coefficients (plh functions) \cite{baker68}, \cite{beynon74}, \cite{beynon77}.

The exponential function (say in base $e$) is an order isomorphism
$$
\exp:M\to((-\infty,0],+,\dotminus,0)
$$
between $M$ and the negative cone of $\Rbb$ endowed with the ordinary sum and the dual truncated difference $b\dotminus a=\min\{0,b-a\}$. 
By~\cite{BlokFerreirim}, $\variety{V}$ is exactly the variety of \newword{cancellative Wajsberg hoops}, i.e., algebras $(A,\tnorm,\to,1)$ satisfying
\begin{gather*}
x\to x = 1 \\
x\tnorm(x\to y) = y\tnorm(y\to x) \\
x\to(y\to z) = (x\tnorm y)\to z \\
x = y\to(x\tnorm y)
\end{gather*}
It can be shown~\cite{BlokFerreirim} that cancellative Wajsberg hoops are categorically equivalent to negative cones of \llgroups.

\begin{theorem}
$\Spec\Free_n(\variety{V})$ is homeomorphic to $\Spec\Free_{n-1}(\variety{MV})$.
\end{theorem}
\begin{proof}
Let $N=(-\infty,0]^n$ be the negative orthant of $\Rbb^n$, and let $P$ be the polyhedral cone spanned positively by $\{b_1e_1+b_2e_2+\cdots+b_{n-1}e_{n-1}+e_n:\vect b1{{n-1}}\in\{0,1\}\}$,
where $(e_1\cdots e_n)$ is the standard basis of $\Rbb^n$.
Let $I$ and $J$ be the principal ideals of $\Fl(n)$ whose elements are all plh functions which are $0$ in $N$ and in $P$, respectively. Form the quotient \llgroups\ $\Fl(n)/I$ and $\Fl(n)/J$, and let $(\Fl(n)/I)^-$, $(\Fl(n)/J)^-$ be their negative cones. $(\Fl(n)/I)^-$ can be identified with the set of all plh functions from $N$ to
$(-\infty,0]$, and analogously for $(\Fl(n)/J)^-$. We have:
\begin{itemize}
\item the categorical equivalence between cancellative Wajsberg hoops and negative cones of \llgroups\ associates $\Free_n(\variety{V})$ to $(\Fl(n)/I)^-$; this is proved, in dual form, in~\cite{cignolitor00}. The spectrum is preserved both by the equivalence and by passing from an \llgroup\ to its negative cone (because every \llgroup\ homomorphism is determined by its behaviour on the negative cone). Therefore $\Spec\Free_n(\variety{V})$ is homeomorphic to $\Spec(\Fl(n)/I)$ (see~\cite[Chapter~10]{bkw} for the spectra of \llgroups);
\item by~\cite[p.~195]{pantiprime}, $\Spec(\Fl(n)/J)$ is homeomorphic to $\Spec\Free_{n-1}(\variety{MV})$.
\end{itemize}
We will establish our claim by showing that $\Fl(n)/I$ and $\Fl(n)/J$ are isomorphic. By~\cite{beynon77} this can be done by triangulating $N$ and $P$ into combinatorially isomorphic complexes of unimodular cones (a polyhedral cone is \newword{unimodular} if it is positively spanned by a $\Zbb$-basis of $\Zbb^n$).

Let $R$ be obtained from the $n\times n$ identity matrix by permuting the first $n-1$ rows. Let $C_R$ be the unimodular cone positively spanned by the columns of
$$
R\begin{pmatrix}
0 & 1 & 1 & \cdots & 1 & 1 \\
0 & 0 & 1 & \cdots & 1 & 1 \\
0 & 0 & 0 & \cdots & 1 & 1 \\
\vdots & \vdots & \vdots & \cdots & \vdots & \vdots \\
0 & 0 & 0 & \cdots & 0 & 1 \\
1 & 1 & 1 & \cdots & 1 & 1
\end{pmatrix}
$$
The collection of all faces of all the $(n-1)!$ $C_R$'s is a complex that triangulates~$P$~\cite[Lemma~2.1]{pantilu}. Analogously, let $D_R$ be the cone spanned by the columns of
$$
R\begin{pmatrix}
0 & -1 & -1 & \cdots & -1 & -1 \\
0 & 0 & -1 & \cdots & -1 & -1 \\
0 & 0 & 0 & \cdots & -1 & -1 \\
\vdots & \vdots & \vdots & \cdots & \vdots & \vdots \\
0 & 0 & 0 & \cdots & 0 & -1 \\
-1 & 0 & 0 & \cdots & 0 & 0
\end{pmatrix}
$$
Again $D_R$ is unimodular and the complex of all $D_R$'s triangulates $N$. This latter fact is easily seen by observing that $D_R$ contains exactly those vectors of $N$ whose first $n-1$ coordinates $\vect\alpha1{{n-1}}$ satisfy $\alpha_{\rho^{-1}1}\le
\alpha_{\rho^{-1}2}\le\cdots\le\alpha_{\rho^{-1}(n-1)}$, where $\rho$ is the permutation that originated $R$. The unimodular complexes thus obtained are clearly combinatorially isomorphic, and this concludes the proof.
\end{proof}

Is is convenient to realize $\Free_n(\variety{V})$ as $(\Fl(n)/I)^-$. Let $\Delta^{n-1}$ be the $(n-1)$-dimensional simplex $\{\sum \alpha_ie_i:\sum \alpha_i=-1\}$. The mapping $\rho$ that associates to a point $u\in\Delta^{n-1}$ the ideal $\rho(u)$ of all plh functions which are $0$ in $u$ is a homeomorphic embedding of $\Delta^{n-1}$ into $\Spec\Free_n(\variety{V})$, whose range is precisely the set of maximal ideals. Let $\sigma:\Free_n(\variety{V})\to\Free_n(\variety{V})$ be a substitution. Then we have a commuting diagram
$$
\begin{CD}
\Delta^{n-1}  @>{\tilde{s}}>>  \Delta^{n-1}  \\
@V{\rho}VV             @VV{\rho}V   \\
\Spec\Free_n(\variety{V})   @>>S>   \Spec\Free_n(\variety{V})
\end{CD}
$$
where $\tilde{s}$ is the piecewise-fractional transformation defined as follows: if $\sigma(x_i)=s_i(\vect x0{{n-1}})$ ($s$ a polynomial in the language $(+,\dotminus,0)$) and $u=(\vect\alpha0{{n-1}})$, then $\tilde{s}(u)=-(\sum s_i(u))^{-1}(s_0(u),\ldots,s_{n-1}(u))$.

On $\Spec\Free_n(\variety{V})$ we have, a priori, two reasonable measures: our standard $\lambda$ obtained by pushing forward via $\pi$ the Lebesgue measure on $\oi^n$, and an unnamed measure obtained by pushing forward via $\rho$ the Lebesgue measure on $\Delta^{n-1}$.

\begin{theorem}
The\label{ref15} unnamed measure coincides with $\lambda$.
\end{theorem}
\begin{proof}
It suffices to prove the statement for the basic closed sets $F_t=\Spec\Free_n(\variety{V})\setminus O_t$. Let us write $t_1$ for the polynomial $t$ written in the language $(\tnorm,\to,1)$, and $t_2$ for the same polynomial written in the language $(+,\dotminus,0)$. Then $\lambda(F_t)$ is the $n$-dimensional Lebesgue measure of $Z(t_1)\subseteq\oi^n$, while the unnamed measure of $F_t$ is the $(n-1)$-dimensional Lebesgue measure, say $r$, of the section $Z(t_2)\cap\Delta^{n-1}$ of the polyhedral cone $Z(t_2)$. The componentwise exponential function gives a diffeomorphism
$$
\exp:(-\infty,0]^n\to(0,1]^n,
$$
and $\exp^{-1}[Z(t_1)]=Z(t_2)$. We have therefore
$$
\lambda(Z(t_1))=\int_{Z(t_2)}\abs{J}\,dx_0\ldots dx_{n-1},
$$
where the Jacobian $J$ of the diffeomorphism has determinant $\exp(x_0+\cdots+x_{n-1})$ in the point $(\vect x0{{n-1}})$. By parametrizing $Z(t_2)$ along a ray, the Riemann integral to the right reduces to the improper integral
$$
\int^0_{-\infty} -x\,r\,\exp(x)\,dx,
$$
which has value $r$; this concludes the proof.
\end{proof}

Theorem~\ref{ref15} makes it possible to formulate the analogue of Theorem~\ref{ref9}: the proof carries over with straightforward modifications.

\begin{theorem}
If\label{ref17} any of the mappings $\tilde{s}$ and $S$ is nonsingular, then so is the other. If this happens, then the systems $(\Delta^{n-1},\tilde{s},\text{Lebesgue measure})$ and $(\Spec\Free_n(\variety{V}),S,\lambda)$ are isomorphic under $\rho$.
\end{theorem}

We conclude our paper with an example. Choose integers $a,b\ge1$, and let $\sigma$ be the substitution over two variables defined by
\begin{equation*}
\begin{split}
x_0 &\mapsto \bigl[\bigl((x_0\to x_1)\to x_1\bigr)\to x_0\tnorm(x_0\to x_1)\bigr]^a \\
x_1 &\mapsto \bigl[(x_0\to x_1)\to x_1\bigr]^b
\end{split}
\end{equation*}
($t^a$ means $t\tnorm t\tnorm\cdots\tnorm t$ $a$ times).
Let $\bar{s}:\oi^2\to\oi^2$ be as in Section~\ref{ref8}, and let $\tilde{s}:\Delta^1\to\Delta^1$ be as in this section. Note that
\begin{equation*}
\begin{split}
s_0 &= \sigma(x_0) = a\bigl[\bigl(x_0+(x_1\dotminus x_0)\bigr)\dotminus\bigl(x_1\dotminus(x_1\dotminus x_0)\bigl)\bigr] \\
&=a\bigl[(x_0\land x_1)\dotminus(x_0\lor x_1)\bigr] = a\bigl[(x_0\land x_1)-(x_0\lor x_1)\bigr] \\
s_1 &= \sigma(x_1) = b\bigl[x_1\dotminus(x_1\dotminus x_0)\bigr] = b\bigl[x_0\lor x_1\bigr]
\end{split}
\end{equation*}
as plh functions from $(-\infty,0]^2$ to $(-\infty,0]$.
It is convenient to identify $\Delta^1$ with $\oi$ via $x\mapsto(-1+x,-x)$. Under this identification $\tilde{s}$ is a transformation on  $\oi$ which depends only on the ratio $q=a/b>0$. Indeed, a straightforward computation shows that
$$
\tilde{s}(x)=\begin{cases}
x\cdot\bigl((1-2q)x+q\bigr)^{-1} & \text{if $0\le x\le 1/2$;} \\
(1-x)\cdot\bigl((1-2q)(1-x)+q\bigr)^{-1} & \text{if $1/2 < x\le 1$.}
\end{cases}
$$
The graph of $\tilde{s}$ is tent-like; since $\tilde{s}(x)=\tilde{s}(1-x)$, it is symmetric with respect to the line $x=1/2$. For $q=1/2$ the slopes of the tent are straight lines, for $1/2<q$ they are convex functions, and for $q<1/2$ are concave. 

\begin{theorem}
Let\label{ref16} $a,b,q,\sigma,\bar{s},\tilde{s}$ be as above.
\begin{enumerate}
\item If $q>1$, then $\lambda$-all points of $\oi$ are attracted to $0$ under iteration of $\tilde{s}$.
\item If $q<1$, then $\tilde{s}$ is exact with respect to a uniquely determined measure~$\ll\lambda$.
\item If $q=1$, then
\begin{itemize}
\item[(3.1)] $\sigma$ is generic;
\item[(3.2)] no point of $(0,1]^2$ has a dense orbit under $\bar{s}$;
\item[(3.3)] $\tilde{s}$ is not ergodic with respect to any probability measure~$\gg\lambda$.
\end{itemize}
\end{enumerate}
\end{theorem}
\begin{proof}
Let $f=\tilde{s}\restriction[0,1/2]$, and note that $\tilde{s}(x)=f(\min(x,1-x))$. We have
$f'(x)=q((1-2q)x+q)^{-2}$; in particular $\abs{f'(0)}=1/q$, and (1) is immediate. We prove (2) by applying Theorem~\ref{ref22}; the only condition that needs checking is the second.
If $1/2\le q$, then $f'(x)\ge f'(0)>1$ is easily shown. Otherwise, if $q<1/2$, then $(1-2q)x+q$ is positive strictly increasing on $[0,1/2]$ and hence $f'(x)$ is positive strictly decreasing on the same interval. Let $d=f^{-1}(1/2)=q/(2q+1)$. We have $\tilde{s}^2(x)=f^2(x)$ on $[0,d]$ and $\tilde{s}^2(x)=f(1-f(x))$ on $[d,1/2]$. By the chain rule, $\abs{(\tilde{s}^2)'}\ge f'(1/2)\cdot f'(d)=(2q+1)^2>1$ on $[0,d]$. We claim that the same bound holds on $[d,1/2]$ (and hence on all $\oi$). By direct computation
\begin{equation*}
\Bigl\lvert\frac{d}{dx}f(1-f(x))\Bigr\rvert=
\frac{1}{\bigl((2q-1)x-(q-1)\bigr)^2},
\end{equation*}
and the denominator describes an upward parabola having vertex in a point $>1/2$. Hence $\abs{(\tilde{s}^2)'}$ is bounded from below in the interval $[d,1/2]$ by the value $(2q+1)^2$ that the displayed expression assumes in $d$, as claimed.

We prove (3) starting from (3.2). Let $p=(\alpha,\beta)\in(0,1]^2$ be on the hyperbola $x_0x_1=c$. Then $\bar{s}(p)=\bigl((\alpha\land\beta)/(\alpha\lor\beta),
\alpha\lor\beta\bigr)$ is on $x_0x_1=d$, where $d=\alpha\land\beta\ge\alpha\beta=c$. It follows that the orbit of $p$ is all above the hyperbola $x_0x_1=c$, and hence it cannot be dense in $(0,1]^2$.

The proof of (3.1) requires a few basic facts about continued fractions: see~\cite{hardywri85} or \cite{Rosen00}. Every real number $\alpha$ has a unique expansion into a continued fraction
$$
\alpha=a_0+\cfrac{1}{a_1+
            \cfrac{1}{a_2+\dotsb
            }}=[a_0,a_1,a_2,\ldots];
$$
the expansion is finite iff $\alpha$ is rational. The \newword{Gauss map} $g:\oi\to\oi$ is defined by $g(0)=0$ and $g(\alpha)=1/\alpha-\lfloor 1/\alpha \rfloor$ if $\alpha\not=0$.
If $\alpha$ has the above expansion, then $g(\alpha)$ has the expansion $[0,a_2,a_3,\ldots]$. Suppose that $\alpha\in\{0\}\cup(1/2,1]$; this means that either $\alpha=[a_0]$ (i.e., $\alpha$ equals $0$ or $1$), or $a_0=0$ and $a_1=1$. Then $\tilde{s}(\alpha)=g(\alpha)$ (by definition if $\alpha\in\{0,1\}$, and because $\tilde{s}(\alpha)=(1-\alpha)\alpha^{-1}=\alpha^{-1}-1=
[a_1,a_2,a_3,\ldots]-1=[0,a_2,a_3,\ldots]=g(\alpha)$ if $\alpha\in(1/2,1)$). Suppose on the other hand $\alpha=[0,a_1,a_2,\ldots]\in(0,1/2]$, so that $a_1\ge2$. Then $\tilde{s}(\alpha)=[0,a_1-1,a_2,\ldots]$; indeed, $(\tilde{s}(\alpha))^{-1}=\alpha^{-1}-1=[a_1,a_2,a_3,\ldots]-1=
[a_1-1,a_2,a_3,\ldots]=[0,a_1-1,a_2,\ldots]^{-1}$.

We see therefore that $\tilde{s}$ acts like a slow continued fraction algorithm. Starting from $\alpha$ and applying $\tilde{s}$ we touch every point which is touched by $g$, but at a slower rate: when we reach, say, $\beta=[0,b_1,b_2,\ldots]$ with $b_1$ a large number, we must move to all points $[0,b_1-c,b_2,\ldots]$ for $1\le c<b_1$, before reaching $g(\beta)$. In particular, for every $\alpha$, the $\tilde{s}$-orbit of $\alpha$ contains the $g$-orbit of $\alpha$. Now, the Gauss map is ergodic with respect to the Gauss measure $\gamma$, which has the same nullsets as the Lebesgue measure~$\lambda$~\cite[Section~4]{Billingsley65}. This implies that $\lambda$-all $\alpha$'s have a dense $g$-orbit, and therefore a dense $\tilde{s}$-orbit as well; by Theorem~\ref{ref17} $\sigma$ is generic.

A proof of (3.3) can be adapted from the proof in~\cite[Remark~6.2.1]{LasotaMac94} for a similar transformation. However, this requires a rather delicate analysis of the convergence of the sequence of densities $\{P^n1\}$, where $P$ is the Perron-Frobenius operator induced by $\tilde{s}$. We give a proof based on continued fractions.
Assume by contradiction that $\tilde{s}$ is ergodic w.r.t.~the probability measure $\mu\gg\lambda$. Let $r=\mu((1/2,1])>0$. By the Birkhoff Ergodic Theorem~\cite[\S1.6]{Walters82}, for $\mu$-all points $\alpha$ we have
\begin{equation}\tag{$*$}
\lim_{n\to\infty}\frac{1}{n}\card\{i:0\le i<n\And
\tilde{s}^i(\alpha)\in[0,1/2]\}=1-r.
\end{equation}
By~\cite[p.~45]{Billingsley65}, for $\lambda$-all points $\alpha=[0,a_1,a_2,\ldots]$ we have
\begin{equation}\tag{$**$}
\lim_{k\to\infty}\frac{k}{a_1+\cdots+a_k}=0.
\end{equation}
Pick an $\alpha$ that satisfies both ($*$) and ($**$), and write $n_k=a_1+\cdots+a_k$; note that the expansion of $\alpha$ must be infinite. Passing to a subsequence in ($*$) we get
\begin{equation*}
\lim_{k\to\infty}\frac{1}{n_k}\card\{i:0\le i<n_k\And
\tilde{s}^i(\alpha)\in[0,1/2]\}=1-r.
\end{equation*}
Since $[0,b_1,b_2,\ldots]\in[0,1/2]$ iff $b_1\ge 2$, we have
$$
\card\{i:0\le i<n_k\And
\tilde{s}^i(\alpha)\in[0,1/2]\}=
\sum_{i=1}^k(a_i-1)=
n_k-k.
$$
We conclude that
\begin{equation*}
\begin{split}
1-r &= \lim_{k\to\infty}(1-k/n_k)=1.
\end{split}
\end{equation*}
This is a contradiction since $r>0$.
\end{proof}

\newcommand{\noopsort}[1]{}

\end{document}